\documentclass[12pt]{amsart}
\usepackage{amssymb,amsmath,latexsym, enumerate, amscd}

\usepackage{color}

\numberwithin{equation}{section}
\theoremstyle{plain}
\newtheorem{thm}{Theorem}[section]

\newtheorem{lem}[thm]{Lemma}
\newtheorem{rmk}[thm]{Remark}
\newtheorem{ex}[thm]{Example}
\newtheorem{cor}[thm]{Corollary}
\newtheorem{prop}[thm]{Proposition}

\newtheorem{df}[thm]{Definition}

\begin{document}
\newcommand{\hgt}{\operatorname{ht}}
\newcommand{\Hht}{\operatorname{Hht}}
\newcommand{\Ext}{\operatorname{Ext}}
\newcommand{\Hom}{\operatorname{Hom}}
\newcommand{\Supp}{\operatorname{Supp}_R}
\newcommand{\Ass}{\operatorname{wAss}_R}
\newcommand{\Att}{\operatorname{Att}_R}
\newcommand{\Min}{\operatorname{Min}_R}
\newcommand{\depth}{\operatorname{depth}}
\newcommand{\pd}{\operatorname{pd}}
\newcommand{\ann}{\operatorname{ann}}
\newcommand{\grade}{\operatorname{grade}}
\newcommand{\Grade}{\operatorname{Grade}}
\newcommand{\pgrade}{\operatorname{p-grade}}
\newcommand{\Spec}{\operatorname{Spec}}
\newcommand{\pdepth}{\operatorname{p-depth}}

\title[Non-Noetherian CM Rings]{Non-Noetherian Cohen-Macaulay Rings}
\author{Tracy Dawn Hamilton}
\address{Department of Mathematics and Statistics\\
California State University, Sacramento\\
Sacramento, CA 95819-6051} \email{hamilton@csus.edu}
\urladdr{http://www.csus.edu/indiv/h/hamiltont}
\author{Thomas Marley}
\address{Department of Mathematics\\
University of Nebraska-Lincoln\\
Lincoln, NE 68588-0130} \email{tmarley1@math.unl.edu}
\urladdr{http://www.math.unl.edu/\textasciitilde tmarley1}

\thanks{The first author was supported for eight weeks during the summer of 2005
through the University of Nebraska-Lincoln's \it Mentoring through
Critical Transition Points \rm grant (DMS-0354281) from the National
Science Foundation.} 
\date{August 5, 2006}

\begin{abstract}  In this paper we investigate a property for commutative rings with identity
which is possessed by every coherent regular ring and is
equivalent to Cohen-Macaulay for Noetherian rings. We study the
behavior of this property in the context of ring extensions (of
various types) and rings of invariants.
\end{abstract}

\maketitle

\section{Introduction}

Over the past several decades Cohen-Macaulay rings have played a
central role in the solutions to many important problems in
commutative algebra and algebraic geometry.  Hochster and Huneke
\cite{HH1} write that for many theorems ``the Cohen-Macaulay
condition (possibly on the local rings of a variety) is just what is
needed to make the theory work.'' However, the study of this
condition has mostly been restricted to the class of (commutative)
Noetherian rings. (Of course, a non-Noetherian ring may be a
Cohen-Macaulay module with respect to some Noetherian subring, but
that is a separate -- though not unrelated-- issue.) The question
has been raised by Glaz (\cite{G2}, \cite{G4}) as to whether there
exists a generalization of the Cohen-Macaulay property to
non-Noetherian rings which has certain desirable properties. In
particular, since many applications of Cohen-Macaulay rings use in
an essential way that regular rings are Cohen-Macaulay, it is
natural to look for a definition of Cohen-Macaulay for arbitrary
rings which has the property that all (Noetherian and
non-Noetherian) regular rings are Cohen-Macaulay. In \cite{G4} Glaz
asks whether such a definition exists, at least for coherent rings.
We give an affirmative answer to this question.

The characterization of Noetherian Cohen-Macaulayness which we
will extend is the following: A Noetherian ring $R$ is
Cohen-Macaulay if and only if
 any sequence $x_1,\dots,x_n$ of elements of $R$ generating a height $n$ ideal
 is a regular sequence.  This characterization will need some modification in the
non-Noetherian case in order to have the desired properties.  In
particular, as heights of ideals can behave erratically in
non-Noetherian rings, we replace the height condition on the
sequence $\bold x=x_1,\dots,x_n$ by conditions on the \u Cech
cohomology $H^n_{\bold x}(R)$ and Koszul homologies $H^i(\bold x;R)$
for $i=1,\dots,n$.  The conditions on the Koszul homologies --
vacuous in the Noetherian case -- ensure that Cech cohomology is
isomorphic to local cohomology. The condition on \u Cech cohomology,
namely that $H^n_{\bold x}(R)_p\neq 0$ for all primes $p$ containing
$(\bold x)$, is equivalent to the condition $\hgt (\bold x)=n$ in
the Noetherian case but can be stronger if the ring is
non-Noetherian. Sequences satisfying these conditions will be called
{\it parameter sequences}.  A parameter sequence such that every
truncation (on the right) by any number of elements is also a
parameter sequence is called a {\it strong parameter sequence}. We
then call a ring {\it Cohen-Macaulay} if every strong parameter
sequence is a regular sequence.

In section four of this paper we show that all coherent regular
rings are Cohen-Macaulay under this definition (Theorem
\ref{regular}). We also give some results concerning the passage of
the Cohen-Macaulay property along ring homomorphisms of special
types (e.g., faithfully flat extensions, localizations, and
quotients by a regular sequence). While we are not able to give
complete answers in all of these situations, it is clear that this
notion of Cohen-Macaulay is less fluid for general rings than for
Noetherian ones.  For instance, we show that a quotient of a
Cohen-Macaulay ring by a non-zero-divisor need not be
Cohen-Macaulay.  On the other hand, we prove a couple of results
which demonstrate the utility of our definition. First, we show that
if $R$ is an excellent Noetherian domain of characteristic $p>0$
then the integral closure $R^+$ of $R$ in an algebraic closure of
the fraction field of $R$ is Cohen-Macaulay (Theorem \ref{Rplus}),
using the difficult result of \cite{HH2} that $R^+$ is a big
Cohen-Macaulay $R$-algebra. Secondly, we show that certain rings of
invariants of coherent regular rings of dimension two are
Cohen-Macaulay (cf. Corollary \ref{invariant}). This is a step
toward resolving a conjecture posed by Glaz in \cite{G2}.

In section two we summarize the basic properties of \u Cech
cohomology, weakly proregular sequences, and non-Noetherian grade.
We establish a connection between non-Noetherian grade and the
vanishing of \u Cech cohomology which mirrors the situation for
classical grade in Noetherian rings (Proposition \ref{CechGrade}).
In section three we define parameter sequences and establish their
basic properties.

Thoughout this article all rings are assumed to be commutative
with identity and all modules unital.  A ring with a unique
maximal ideal is called `quasi-local' while the term `local' is
reserved for Noetherian quasi-local rings.  We let $\bold x$
denote a finite sequence of elements $x_1,\dots,x_{\ell}$ of ring
$R$. The length of the sequence $\bold x$ is denoted by $\ell
(\bold x)$. Given any sequence $\bold x$ of $R$  we let $\bold x'$
denote the sequence obtained by truncating the last element from
$\bold x$; i.e, $\bold x'= x_1,\dots,x_{\ell-1}$ where
$\ell=\ell(\bold x)$.

\section{\u Cech cohomology and non-Noetherian grade}

Let $R$ be a ring and $x$ an element of $R$. Let $C(x)$ denote the
cochain complex
$$0\rightarrow  R \xrightarrow[r\rightarrow \frac{r}{1}]{} R_x\rightarrow 0$$
where the position of $R$ is in degree zero and the differential is
the natural localization map. For a sequence $\bold x$ of elements
of $R$, the \u Cech complex $C(\bold x)$ is inductively defined by
$C(\bold x):=C(\bold x')\otimes_R C(x_{\ell})$, where
$\ell=\ell(\bold x)$. If $M$ is an $R$-module, then we set $C(\bold
x;M):=C(\bold x)\otimes_RM$. The $i$th {\it \u Cech cohomology}
$H^i_{\bold x}(M)$ of $M$ with respect to the sequence $\bold x$ is
defined to be the $i$th cohomology of $C(\bold x;M)$.

We list here some of the elementary properties of \u Cech
cohomology. Proofs of these results are either elementary or can
be found in section 5.1 of \cite{BS}.

\begin{prop} \label{Cech-prop} Let $R$ be a ring, $\bold x$ a finite sequence
of elements of $R$ and $M$ an $R$-module.
\begin{enumerate}[(a)]
\item $H^i_{\bold x}(M)=0$ for $i<0$ and $i>\ell(\bold x)$. \item
Given a short exact sequence $0\to A\to B\to C\to 0$ of
$R$-modules there exists a natural long exact sequence $$\cdots
\to H^i_{\bold x}(A)\to H^i_{\bold x}(B)\to H^i_{\bold x}(C)\to
H^{i+1}_{\bold x}(A)\to \cdots. $$\item There exists a long exact
sequence
$$
\cdots \to H^i_{\bold x}(M) \to H^i_{\bold x'}(M)\xrightarrow{f_i}
H^i_{\bold x'}(M)_{x_{\ell}}\to H^{i+1}_{\bold x}(M)\to \cdots $$
where $f_i$ is (up to a sign) the natural localization map. \item
For all $i$ the module $H^i_{\bold x}(M)$ is $(\bold x)R$-torsion;
i.e., every element is annihilated by a power of $(\bold x)R$.
 \item
If $\bold y$ is a finite sequence of elements of $R$ such that
$\sqrt{(\bold x)R}=\sqrt{(\bold y)R}$, then $H^i_{\bold x}(M)\cong
H^i_{\bold y}(M)$ for all $i$. \item (Change of Rings) Let $f:R\to
S$ be a ring homomorphism and $N$ an $S$-module. Then $H^i_{\bold
x}(N)\cong H^i_{f(\bold x)}(N)$ for all $i$. \item (Flat Base
Change) Let $f:R\to S$ be a flat ring homomorphism and $M$ an
$R$-module. Then $H^i_{\bold x}(M)\otimes_RS\cong H^i_{\bold
x}(M\otimes_RS)\cong H^i_{f(\bold x)}(M\otimes_RS)$ for all $i$.
\item $H^{\ell(\bold x)}_{\bold x}(M)\cong H^{\ell(\bold
x)}_{\bold x}(R)\otimes_RM$ and $\Supp H^{\ell(\bold x)}_{\bold
x}(M)\subseteq \Supp M/(\bold x)M$.
\end{enumerate}
\end{prop}

For an ideal $I$ of $R$ the $i$th {\it local cohomology} $H^i_I(M)$
of $M$ with support in $V(I)$ is defined by
$$H^i_I(M):=\lim_{\longrightarrow}\Ext^i_R(R/I^n,M).$$  In the case $R$ is Noetherian
we have $H^i_{\bold x}(M)\cong H^i_{I}(M)$ for all $i$, where
$I=(\bold x)R$.   However, local cohomology and \u Cech cohomology
are not in general isomorphic over non-Noetherian rings. Schenzel
\cite{Sch} gives necessary and sufficient conditions on a sequence
$\bold x$ such that $H^i_{\bold x}(M)\cong H^i_{I}(M)$ for all $i$
and $R$-modules $M$:  For $x\in R$ let $K(x)$ denote the Koszul
chain complex $0\to R\xrightarrow{x}R\to 0$, where the first $R$ is
in degree $1$.  For a sequence $\bold x=x_1,\dots,x_{\ell}$ the
Koszul complex $K(\bold x)$ is defined to be the chain complex
$K(x_1)\otimes \cdots \otimes K(x_{\ell})$.  We denote the homology
of $K(\bold x)$ by $H_i(\bold x)$.  For $m\ge n$ there exists a
chain map $\phi^m_n(\bold x):K(\bold x^m)\to K(\bold x^n)$ given by
$\phi^m_n(\bold x)=\phi^m_n(x_1)\otimes \cdots \otimes
\phi^m_n(x_{\ell})$, where $\phi^m_n(x)$ is the chain map
$$
\begin{CD}
0 @>>> R @>x^m>> R @>>> 0 \\
@. @VV x^{m-n}V @VV=V @. \\
0 @>>> R @>x^n>> R @>>> 0.
\end{CD}
$$
Hence, $\{K(\bold x^m),\phi^m_n(\bold x)\}$ is an inverse system of
chain complexes.  The sequence $\bold x$ is called {\it weakly
proregular} if for each $n$ there exists an $m\ge n$ such that the
maps $\phi^m_n(\bold x)_*:H_i(\bold x^m)\to H_i(\bold x^n)$ are zero
for all $i\ge 1$ \cite{Sch}.  Note that an element $x$ is weakly
proregular if and only if there exists an $n\ge 1$ such that
$(0:x^n)=(0:x^{n+1})$.  We note the following elementary remarks:

\begin{rmk} \label{proregular} Let $R$ be a ring and $\bold x$ a finite sequence of
elements from $R$.
\begin{enumerate}[(a)]
\item If $\bold x$ is weakly proregular then so is any
permutation of $\bold x$. \item Any regular sequence is proregular.
\item Suppose $f:R\to S$ is a flat ring homomorphism.  If $\bold x$
is weakly proregular on $R$ then $f(\bold x)$ is weakly proregular
on $S$.  The converse is true if $f$ is faithfully flat.
\end{enumerate}
\end{rmk}

 We will need the following result, due to Schenzel:

\begin{thm} \label{Schenzel} \cite[Theorem 3.2]{Sch}  Let $R$ be a commutative ring, $\bold
x$ a finite sequence of elements of $R$, and $I=(\bold x)R$. The
following conditions are equivalent:
\begin{enumerate}[(a)]
\item $\bold x$ is weakly proregular. \item $H^i_{\bold x}(E)=0$
for all injective $R$-modules $E$ and $i\neq 0$. \item For every
$R$-module $M$ there exist natural isomorphisms $H^i_{\bold
x}(M)\to H^i_I(M)$ for all $i$.
\end{enumerate}
\end{thm}

As a consequence of this theorem we obtain that any finite sequence
of elements in a Noetherian ring is weakly proregular. Also, if
$\sqrt{(\bold x)R}=\sqrt{(\bold y)R}$ where $\bold x$ and $\bold y$
are finite sequences (but not necessarily of the same length), then
$\bold x$ is weakly proregular if and only if $\bold y$ is.  As an
example of an element in a ring which is not weakly proregular,
consider the image of $x$ in the ring $S=\mathbb
Z[x,y_1,y_2,y_3,\dots]/(xy_1,x^2y_2,x^3y_3,\dots)$. By Schenzel's
theorem, there exists an injective $S$-module $E$ such that
$H^1_x(E)\not\cong H^0_{(x)}(E)=0$.

Let $R$ be a ring of Krull dimension $d$ and $\bold x$ a sequence
of elements of $R$.  Two central results concerning \u Cech
cohomology over a Noetherian ring are:
\begin{enumerate}
\item $H^i_{\bold x}(M)=0$ for $i>d$. \item If $R$ is local and
$\dim R/(\bold x)=0$ then $H^d_{\bold x}(R)\neq 0$.
\end{enumerate}
Part (2) may fail for non-Noetherian rings.  For example, suppose
$V$ is a  valuation domain of finite dimension $d>1$.  As the prime
ideals of $V$ are linearly ordered, there exists a non-unit $x\in V$
such that $\dim V/xV=0$, while clearly $H^d_x(V)=0$.  On the other
hand, (1) holds for arbitrary commutative rings:

\begin{prop} \label{vanishing} Let $R$ be a ring of finite dimension $d$ and
$\bold x$ a sequence of elements from $R$.  Then $H^i_{\bold
x}(M)=0$ for all $i>d$.
\end{prop}

{\it Proof:} We use induction on $d$.  It is enough to prove the
statement in the case $(R,m)$ is a quasi-local ring and $(\bold
x)R\subseteq m$. If $d=0$ then every element of $m$ is nilpotent.
Hence, $H^i_{\bold x}(M)=0$ for all $i\ge 1$.  Now suppose that
$d>0$ and that the proposition holds for all rings of dimension at
most $d-1$. Let $\ell=\ell(\bold x)$. For $\ell\le d$ there is
nothing to show. Suppose $\ell>d$ and let $j$ be the largest
integer such that $H^j_{\bold x}(M)\neq 0$. Assume $j>d$. By
induction on $\ell$, $H^j_{\bold x'}(M)=0$. From the exact
sequence
$$
\dots \to H^{j-1}_{\bold x'}(M) \to H^{j-1}_{\bold
x'}(M)_{x_{\ell}} \to H^j_{\bold x}(M) \to 0
$$
we obtain that $H^{j-1}_{\bold x'}(M_{x_{\ell}})\neq 0$.  As $R$
is quasi-local of dimension $d$, $\dim R_{x_{\ell}}\le d-1$.
Thus, $j-1\le d-1$, a contradiction. \qed

\medskip

We now briefly discuss non-Noetherian grade, also referred to as
`polynomial grade' \cite{EN} and `true grade' \cite{No}. This
notion dates back to the early 1970s (e.g., \cite{Ba}, \cite{Ho})
in connection with the study of finite free resolutions over
(arbitrary) commutative rings. Hochster appears to have been the
first to notice that the pathological behavior of `classical'
grade in the non-Noetherian case can be remedied by adjoining
indeterminates to the ring (\cite[footnote, p.132]{No}).

To begin we recall some terminology from \cite{Ho}.  Let $R$ be a
ring and $M$ an $R$-module.  A sequence $\bold x=x_1,\dots
x_{\ell}\in R$ is called a {\it possibly improper regular
sequence} on $M$ if $x_i$ is a non-zero-divisor on
$M/(x_1,\dots,x_{i-1})M$ for $i=1,\dots,\ell$. If in addition
$M\neq (\bold x)M$ we call $\bold x$ a {\it regular sequence} on
$M$. Given an ideal $I$ of $R$, the {\it classical grade} of $I$
on $M$, denoted $\grade (I,M)$, is defined to be the supremum of
the lengths of all possibly improper regular sequences on $M$
contained in $I$. In the case $R$ is Noetherian and $M$ is
finitely generated, $\grade (I,M)>0$ if and only if $(0:_MI)=0$.
However, there are examples of finitely generated ideals in
non-Noetherian rings which have annihilator zero but consist
entirely of zero-divisors on the ring (cf. \cite{V} or Example
\ref{trivial-pgrade} below). This phenomenon disappears if one
first passes to a polynomial ring extension of $R$. The following
lemma is the central insight behind polynomial grade:

\begin{lem} \label{Hochster} Let $R$ be a ring, $I=(x_1,\dots,x_{\ell})R$, and
$M$ an $R$-module.  Then $\grade (IR[t],R[t]\otimes_RM)>0$ if and
only if $(0:_MI)=0$. In particular, $(0:_MI)=0$ if and only if
$x_1+x_2t+\cdots +x_{\ell}t^{{\ell-1}}\in IR[t]$ is a
non-zero-divisor on $R[t]\otimes_RM$.
\end{lem}

{\it Proof:} See Chapter 5, Theorem 7 of \cite{No}. \qed

\medskip
For an ideal $I$ of $R$ and $R$-module $M$, the {\it polynomial
grade} of $I$ on $M$ is defined by
$$\pgrade (I,M):=\lim_{m\to \infty} \grade (IR[t_1,\dots,t_m],R[t_1,\dots,t_m]\otimes_RM).$$
It is easily seen (cf. \cite{A}, \cite{Ho}) that
$$\pgrade (I,M)=\sup \{\grade
(IS,S\otimes_RM)\mid S \ \text{a faithfully flat $R$-algebra}
\}.$$ We denote $\grade(I,R)$ and $\pgrade (I,R)$ by $\grade I$
and $\pgrade I$, respectively.  If $(R,m)$ is quasi-local, then
$\grade (m,M)$ and $\pgrade (m,M)$ are denoted by $\depth M$ and
$\pdepth M$, respectively. We note that if $R$ is Noetherian then
$\grade (I,M)=\pgrade (I,M)$ for all ideals $I$ and finitely
generated $R$-modules $M$. The following proposition summarizes
the essential properties of polynomial grade.  Proofs of these
results can be found in Chapter 5 of \cite{No}.

\newpage

\begin{prop} \label{polygrade}  Let $R$ be a ring and $I$ an ideal of $R$.
\begin{enumerate}[(a)]
\item If $\bold x$ is a regular sequence on $M$ contained in $I$
then
$$\pgrade (I,M)=\pgrade(I,M/(\bold x)M)+\ell(\bold x).$$ \item $\pgrade (I,M)= \pgrade
(P,M)$ for some prime ideal $P$ containing $I$.  In particular,
$\pgrade(I,M)=\pgrade (\sqrt{I},M)$. \item $\pgrade (I,M)=\sup
\{\pgrade (J,M) \mid J\subseteq I,\ J \text{\rm \ a finitely
generated ideal}\}$. \item If $I$ is generated by $n$ elements and
$IM\neq M$, then $\pgrade I\le n$. Furthermore, $\pgrade
(I,M)=\grade (IR[t_1,\dots,t_n],R[t_1,\dots,t_n]\otimes_RM)$.
\end{enumerate}
\end{prop}

For a sequence $\bold x$ of $R$ and an $R$-module $M$, the $i$th
Koszul homology of $\bold x$ on $M$, denoted $H_i(\bold x;M)$, is
defined to be the $i$th homology of $K(\bold x)\otimes_RM$.  The
following proposition relates polynomial grade with the vanishing of
Koszul homology and \u Cech cohomology:

\begin{prop} \label{CechGrade} Let $R$ be a ring, $\bold x$ a
finite sequence of elements from $R$ of length $\ell=\ell(\bold
x)$, $I=(\bold x)R$, and $M$ an $R$-module.  The following
integers (including the possibility of $\infty$) are equal:
\begin{enumerate}
\item $\pgrade (I,M);$ \item $\sup \{k\ge 0\mid H_{\ell-k}(\bold
x,M)=0\ \text{ for all } i<k\};$ \item $\sup \{k\ge 0\mid
H^i_{\bold x}(M)=0 \text{ for all } i<k\}.$
\end{enumerate}
Moreover, $IM\neq M$ if and only if any one of the above integers
is finite.
\end{prop}

{\it Proof:} The equality of (1) and (2) was established by Barger
\cite{Ba} and Alfonsi \cite{A}.  We prove the equality of (1) and
(3).  Let $p=\pgrade (I,M)$ and $h=h(\bold x,M)$ the quantity
representing (3). We first assume $p<\infty$ and use induction on
$p$ to prove $p=h$. If $p=0$ then by Lemma \ref{Hochster} we have
$(0:_MI)\neq 0$. Hence, $H^0_{\bold x}(M)\neq 0$ and $h=0$. Assume
now that $p>0$. Then $\grade (IR[t],R[t]\otimes_RM)>0$. As $R\to
R[t]$ is a faithfully flat ring extension, it follows by
Proposition \ref{Cech-prop} that $H^i_{\bold x}(M)=0$ if and only
if $H^i_{\bold x}(R[t]\otimes_RM)=0$ for all $i$. Hence, by
replacing $R$ by $R[t]$, we may assume $\grade (I,M)>0$. Let $u\in
I$ be a non-zero-divisor on $M$. Since $\pgrade (I,M/uM)=p-1$, we
have by induction that $h(\bold x,M/uM)=p-1$. The short exact
sequence $0\to M\xrightarrow{u} M\to M/uM\to 0$ induces the long
exact sequence
$$\dots \to H^{i-1}_{\bold x}(M/uM)\to H^i_{\bold x}(M)\xrightarrow{u} H^i_{\bold x}(M) \to \cdots.$$ Thus, for
$i\le p-1$, multiplication by $u$ on $H^i_{\bold x}(M)$ is
injective.  However, since $u\in I$ every element of $H^i_{\bold
x}(M)$ is annihilated by a power of $u$.  This implies $H^i_{\bold
x}(M)=0$ for $i\le p-1$.  Finally, the same long exact sequence
yields the exactness of
$$0\to H^{p-1}_{\bold x}(M/uM)\to H^p_{\bold x}(M),$$
which implies $H^p_{\bold x}(M)\neq 0$.  Thus, $h=p$.

Suppose $h<\infty$.  As before, if $h=0$ then $p=0$.  If $h>0$
then by Lemma \ref{Hochster} we have $\grade
(IR[t],R[t]\otimes_RM)>0$. Since $h(I,M)=h(IR[t],
R[t]\otimes_RM)$, we may assume the existence of $u\in I$ which is
a non-zero-divisor on $M$. Then $\pgrade (I,M/uM)=p-1$ and (using
the same long exact sequence as above) $h(I,M/uM)=h-1$.  By
induction, $h-1=p-1$. \qed \medskip

Let $M$ be an $R$-module.  A prime ideal $P$ is said to be {\it
weakly associated} to $M$ if $P$ is minimal over $(0:_Rx)$ for
some $x\in M$ (cf. \cite{Bk}). We denote the set of weakly
associated primes of $M$ by $\Ass M$.  It is easily seen that if
$R$ is Noetherian $\Ass M= \operatorname{Ass}_RM$ for all
$R$-modules $M$. As in the Noetherian case, the union of the
weakly associated primes of $M$ is the set of zero-divisors on $M$
and $\Ass M=\emptyset$ if and only if $M=0$.  We prove the
following elementary (and presumably well-known) result:

\begin{lem} \label{attached} Let $M$ be an $R$-module and $p\in \Ass M$.  Then
$\pdepth_{R_p}M_p=0$.
\end{lem}

{\it Proof:} By localizing at $p$, we may assume $(R,m)$ is
quasi-local and $m=\sqrt{(0:_Rx)}$ for some $x\in M$. Let $J$ be a
finitely generated ideal contained in $m$.  Then $J^n\subseteq
(0:_Rx)$ for some $n$.  Hence, $(0:_MJ^n)\neq 0$ which implies
$\pgrade (J,M)=\pgrade (J^n,M)=0$ by Lemma \ref{Hochster} and
Proposition \ref{polygrade}(b).  By part (c) of Proposition
\ref{polygrade}, we obtain $\pgrade (m,M)=0$. \qed
\medskip

While polynomial grade has many of the same properties as
classical grade for Noetherian rings, one important difference is
that a ring may contain ideals of polynomial grade $j>1$ but no
ideals of polynomial grade $i$ for $0<i<j$.  To see this, we first
prove the following proposition, which is adapted from \cite{V}:

\begin{prop} \label{Wolmer}  Let $(R,m)$ be a quasi-local ring
of dimension $d$.  Fix an integer $i\ge 0$ and let
\begin{displaymath}
M_i:=\bigoplus_{\substack{p\in \Spec R \\ \hgt p\le i}}k(p)
\end{displaymath} where $k(p)$ is the residue field of $R_p$. Let $S=R\times M_i$ be the trivial extension of
$R$ by $M_i$, $j:S\to R$ the natural projection, and $I$ a
finitely generated ideal of $S$. Then $\hgt I=\hgt j(I)$ and
$$\pgrade I=\left\{ \begin{array}{ll} 0 & \textrm{if
$\hgt I\le i$,} \\
\pgrade j(I) & \textrm{if $\hgt I>i$.}
\end{array} \right.
$$
Moreover, for any sequence $\bold x$ of $S$, $\bold x$ is weakly
proregular on $S$ if and only if $j(\bold x)$ is weakly proregular
on $R$.
\end{prop}

{\it Proof:}  Since the ideal $0\times M$ is nilpotent, there is a
bijective correspondence between $\Spec S$ and $\Spec R$ given by
$P\to j(P)$.   Hence, $\hgt I=\hgt j(I)$ and
$\sqrt{I}=\sqrt{j(I)S}$.  By Proposition \ref{polygrade}(b),
$\pgrade I=\pgrade j(I)S=\pgrade (j(I),S)$.  Thus, it suffices to
show that for all finitely generated ideals $J$ of $R$, $\pgrade
(J,S)=0$ if $\hgt J\le i$ and $\pgrade (J,S)=\pgrade (J,R)$ if
$\hgt J>i$. Suppose first that $\hgt J\le i$.  Then $J$ is
contained in some prime $p$ of $R$ of height $j\le i$.  Let
$\alpha$ be an element of $M_i$ which is nonzero in the component
corresponding to $k(p)$ and zero in all other components. Clearly,
$Js=0$ where $s=(0,\alpha)\in S$.  Hence, $\pgrade (J,S)=0$ by
Lemma \ref{Hochster}.

Suppose that $\hgt J>i$.  Let $J=(\bold x)R$.  For any prime $p$ of
height less than $\hgt J$ we have $(\bold x)k(p)=k(p)$.  By the
change of rings isomorphism (Proposition \ref{Cech-prop}(f)),
$H^i_{\bold x}(k(p))=0$ for all $i\ge 0$.  Therefore, $H^i_{\bold
x}(M_i)=0$ for all $i$.  As $S\cong R\oplus M_i$ as $R$-modules, we
have $H^i_{\bold x}(S)\cong H^i_{\bold x}(R)$ for all $i$. Thus,
$\pgrade (J,S)=\pgrade (J,R)$ by Proposition \ref{CechGrade}.

To prove the last statement, note that as $\sqrt{(\bold
x)S)}=\sqrt{j(\bold x)S}$, we have by Proposition \ref{Schenzel}
that $\bold x$ is weakly proregular on $S$ if and only if $j(\bold
x)$ is weakly proregular on $S$.  Thus, it suffices to prove that if
$\bold x$ is a finite sequence of elements from $R$, then $\bold x$
is weakly proregular on $R$ if and only if it is weakly proregular
on $S$.  However, as $R$-modules $H_j(\bold x^n;S)\cong H_j(\bold
x^n;R)\oplus H_j(\bold x^n;M_i)$ for all$j$.  Since $M_i$ is a
direct sum of fields, it is easy to see that the maps $H_j(\bold
x^{n+1};M_i)\to H_j(\bold x^n; M_i)$ are zero for all $j\ge 1$ and
$n\ge 0$.  (Note that if $F$ is a field and $\bold y$ a sequence in
$F$, then $H_i(\bold y,F)\neq 0$ for some $i$ if and only if $\bold
y$ is the zero sequence.)  Hence, $\bold x$ is weakly proregular on
$R$ if and only if it is weakly proregular on $S$. \qed
\medskip

Applying this Proposition in the case $R$ is a Cohen-Macaulay local
ring, we get the following:

\begin{ex} \label{trivial-pgrade} Let $(R,m)$ be a Cohen-Macaulay local ring of dimension
$d>0$.  Let $S=R\times M_{d-1}$ as in Proposition \ref{Wolmer}.
Then $S$ is a quasi-local ring of dimension $d$ with maximal ideal
$n=m\times M_{d-1}$ with the following properties:
\begin{enumerate}[(a)]
\item $\pdepth S=\dim S$. \item $\pgrade I=0$ for all ideals $I$
of $S$ such that $\sqrt{I}\neq n$; in particular, $n$ consists
entirely of zero-divisors.  \item $\pdepth S_P=0$ for all $P\in
\Spec S\setminus \{n\}$.
\end{enumerate}
\end{ex}

{\it Proof:} As we noted in the proof of Proposition \ref{Wolmer},
$\sqrt{I}=\sqrt{j(I)S}$ for every ideal $I$ of $S$.  As $j(I)$ is
finitely generated (since $R$ is Noetherian), we see that
Proposition \ref{Wolmer} applies to all ideals of $S$.  Parts (a)
and (b) now follow.  For part (c), note that for any $P\in \Spec
S$, $S_P\cong R_{j(P)}\times (M_{d-1})_{j(P)}$.  \qed

\medskip

We end this section with a statement of the Auslander-Buchsbaum
theorem for quasi-local rings.   We denote the projective
dimension of an $R$-module $M$ by $\pd_RM$.  A {\it finite free
resolution} (FFR) of $M$ is a resolution of $M$ of finite length
consisting of finitely generated free modules in each degree.

\begin{prop} \label{ABT} Let $(R,m)$ be a quasi-local ring
and $M$ an $R$-module which has an FFR.  Then
$$\pd_R M + \pdepth M= \pdepth R.$$
\end{prop}

{\it Proof:} See Chapter 6, Theorem 2 of {\cite{No}. \qed

\section{Parameter sequences}

A sequence of elements $\bold x$ in a (Noetherian) local ring $R$
is said to be a system of parameters (s.o.p.) if $\hgt (\bold
x)R=\ell (\bold x)=\dim R$. If $\hgt (\bold x)R=\ell(\bold x)<\dim
R$, we say $\bold x$ is a partial s.o.p.  We wish to extend this
notion to sequences in non-Noetherian rings using homological
properties of the ring instead of height conditions.

\begin{df} {\rm Let $R$ be a ring and $M$ an $R$-module.  A finite sequence $\bold x$
of elements of $R$ is called a {\it parameter sequence} on $R$
provided the following conditions hold:
\begin{enumerate}
\item $\bold x$ is a weakly proregular sequence; \item $(\bold
x)R\neq R$; \item $H^{\ell(\bold x)}_{\bold x}(R)_p\neq 0$ for all
prime ideals $p$ containing $(\bold x)R$.
\end{enumerate}
A parameter sequence of length one on $R$ is called a {\it
parameter} of $R$. The sequence $\bold x$ is called a {\it strong
parameter sequence} on $R$ if $x_1,\dots,x_i$ is a parameter
sequence on $R$ for $i=1,\dots,\ell(\bold x)$.}
\end{df}

We define the ideal generated by the empty sequence to be the zero
ideal. Thus, the empty sequence is a parameter sequence of length
zero on any ring. The empty sequence will also be considered as a
regular sequence of length zero on any ring. These conventions
will allow us to begin proofs by induction with $\ell(\bold x)=0$.

The following remark shows that parameter sequences coincide with
(partial) systems of parameters if the ring is Noetherian.

\begin{rmk} \label{Noeth-remark} Let $R$ be a Noetherian ring and $\bold
x$ a finite sequence of elements of $R$.  Then $\bold x$ is a
parameter sequence on $R$ if and only if $\hgt(\bold x)R=\ell(\bold
x)$.
\end{rmk}

{\it Proof:}  Recall that any sequence of elements in a Noetherian
ring is a weakly proregular sequence.  Also, by convention $\hgt
(\bold x)=\infty$ if and only if $(\bold x)R=R$.  Let $p$ be a prime
of height $h$ which is minimal over $(\bold x)$ and let $\ell =
\ell(\bold x)$. By Krull's Principle Ideal Theorem we have $h\le
\ell$. Since $(\bold x)R_p$ is primary to $pR_p$ and using standard
facts about local cohomology, we get that $H^\ell_{\bold
x}(R)_p\cong H^\ell_{pR_p}(R_p)\neq 0$ if and only if $h=\ell$. \qed

\medskip

We cite some elementary properties of parameter sequences:

\begin{prop} \label{param-prop} Let $R$ be a ring and $\bold x$ a finite sequence of elements of $R$.
\begin{enumerate}[(a)]
\item Any permutation of a parameter sequence on $R$ is again a
parameter sequence on $R$. \item If $\sqrt{(\bold y)R}=\sqrt{(\bold
x)R}$ and $\ell(\bold y)=\ell(\bold x)$ then $\bold x$ is a
parameter sequence on $R$ if and only if $\bold y$ is a parameter
sequence on $R$. \item Let $f:R\to S$ be a flat ring homomorphism.
If $\bold x$ is a (strong) parameter sequence on $R$ and $S/(f(\bold
x))S\neq 0$ then $f(\bold x)$ is a (strong) parameter sequence on
$S$.  The converse holds if $f$ is faithfully flat.  \item Let
$f:R\to S$ be a ring homomorphism and $\bold x$ a weakly proregular
sequence on $R$ such that $(f(\bold x))S_p\neq S_p$ for all primes
$p$ of $R$ minimal over $(\bold x)R$. If $f(\bold x)$ is a parameter
sequence on $S$ then $\bold x$ is a parameter sequence on $R$. \item
If $\pgrade ((\bold x),R)=\ell(\bold x)$ then $\bold x$ is a
parameter sequence on $R$. \item Every regular sequence on $R$ is a
strong parameter sequence.
\end{enumerate}
\end{prop}

{\it Proof:}  Throughout the proof, let $\ell=\ell(\bold x)$. Parts
(a) and (b) follow readily from Proposition \ref{Cech-prop} and
Remark \ref{proregular}. For part (c), we have $f(\bold x)$ is
proregular by Remark \ref{proregular}.  Let $Q$ be a prime of $S$
containing $(f(\bold x))S$. Then $p=f^{-1}(Q)$ is a prime of $R$
containing $(\bold x)R$.  Since the map $R_p\to S_Q$ is faithfully
flat, we obtain
$$H^{\ell}_{f(\bold x)}(S)_Q\cong H^{\ell}_{\bold
x}(R)_p\otimes_{R_p}S_Q\neq 0.$$  The converse is proved
similarly, again using Remark \ref{proregular}.

For part (d), let $p$ be a prime of $R$ minimal over $(\bold x)$.
Since $(f(\bold x))S_p\neq S_p$ and $f(\bold x)$ is a parameter
sequence, we have
$$H^{\ell}_{\bold x}(R)_p\otimes_R S\cong H^{\ell}_{f(\bold x)}(S)_p\neq
0.$$  Hence, $H^{\ell}_{\bold x}(R)_p\neq 0$.

To prove (e), we first note that $H_i(\bold x^n)=0$ for all $i\ge 1$
by Propositions \ref{polygrade}(b) and \ref{CechGrade}.  Hence,
$\bold x$ is weakly proregular.  Next, note that $\pgrade ((\bold
x)R_p,R_p)<\infty$ if and only if $(\bold x)R_p\neq R_p$ by
Proposition \ref{CechGrade}. Since localization does not decrease
p-grade and p-grade (if finite) is bounded above by the length of
the sequence, we see that $\pgrade ((\bold x)R_p)=\ell$ for all
primes $p$ containing $(\bold x)R$. Therefore, $H^{\ell}_{(\bold
x)}(R)_p\neq 0$ for all $p\supseteq (\bold x)R$ by Proposition
\ref{CechGrade}. Part (f) is an immediate consequence of (e). \qed

\medskip

The following lemma allows us to give a simple characterization of
parameters on $R$:

\begin{lem} \label{easylemma} Let $R$ be a ring and $x\in J(R)$, where $J(R)$ is the
Jacobson radical of $R$.  Then $H^1_x(R)=0$ if and only if $x$ is
nilpotent.
\end{lem}

{\it Proof:} Suppose $H^1_x(R)\cong R_x/R=0$.  Then there exists an
$r\in R$ such that $\frac{1}{x}=\frac{r}{1}$ in $R_x$.  Thus for
some $i$, $(1-rx)x^i=0$.  As $x\in J(R)$, $1-rx$ is a unit and hence
$x^i=0$. \qed

\medskip

\begin{cor} \label{parameter} Let $R$ be a ring and $x\in R$ a nonunit.  Then $x$ is
a parameter on R if and only if $\hgt xR\ge 1$ and
$(0:x^n)=(0:x^{n+1})$ for some $n\ge 1$.
\end{cor}

{\it Proof:} Let $p$ be a prime minimal over $xR$.  By Lemma
\ref{easylemma}, $H^1_x(R)_p=0$ if and only if $x$ is nilpotent in
$R_p$, which is the case if and only if $\hgt p=0$.  Hence,
$H^1_x(R)_p\neq 0$ for all primes $p$  minimal over $xR$ if and only
if $x$ is not in any minimal prime of $R$.  As noted in the
paragraph preceding Remark \ref{proregular}, $x$ is weakly
proregular if and only if $(0:x^n)=(0:x^{n+1})$ for some $n$. \qed

\medskip

As a consequence of Proposition \ref{vanishing} we have:

\begin{prop} \label{height} Let $R$ be a ring and $\bold x$ a parameter sequence on
$R$.  Then $\hgt (\bold x)R\ge \ell(\bold x)$. \end{prop}

\medskip

The following example shows that the converse to this proposition
can be false even in the case the sequence in weakly proregular and
$R$ is a coherent regular ring:

\begin{ex} \label{valuation}{\rm Let $V$ be a valuation domain of
(Krull) dimension $2$.  Let $m$ be the maximal ideal of $V$ and $P$
the (unique) prime ideal lying between $(0)$ and $m$.  Choose a
non-zero element $x\in P$ and $y\in m\setminus P$.  Then $x,y$ is
weakly proregular and $\hgt (x,y)V=2$ but $x,y$ is not a parameter
sequence.} \end{ex}

{\it Proof:} Clearly $\hgt (x,y)V=\hgt m=2$.  Since $(x,y)V=yV$,
$H^2_{x,y}(V)=0$.  Thus, $x,y$ is not a parameter sequence. It
suffices to show that $x,y$ is weakly proregular.  As $V$ is a
domain, $H_2(x^n,y^n)=0$ for all $n$. Let $\alpha\in
H_1(x^{2n},y^{2n})$ and $(r,s)\in R^2$ a lifting of $\alpha$. Then
$rx^{2n}+sy^{2n}=0$. The map $H_1(x^{2n},y^{2n})\to H_1(x^n,y^n)$ is
induced by the map $(r,s)\to (rx^n,sy^n)$.  As $x\in yV$, $x=b y$
for some $\beta \in V$. Let $a=rb^n$. Then
\begin{align}
rx^n&=rb^ny^n=a y^n, \notag \\
sy^n&=\frac{sy^{2n}}{y^n}=-\frac{rx^{2n}}{y^n}=-\frac{rb^n y^n
x^n}{y^n}=-a x^n. \notag
\end{align}
Hence, $(rx^n,sy^n)=a(y^n,-x^n)$ is a boundary.  Thus, the map
$H_1(x^{2n},y^{2n})\to H_1(x^n,y^n)$ is zero. \qed

\medskip

We end this section by characterizing strong parameter sequences
on trivial extensions of the type described in Proposition
\ref{Wolmer}.

\begin{prop} \label{sps-trivial} Let the notation be as in the statement of
Proposition \ref{Wolmer}. The following are equivalent for a finite
sequence $\bold x$ of $S$:
\begin{enumerate}[(a)] \item $\bold x$ is a
(strong) parameter sequence on $S$.  \item $j(\bold x)$ is a
(strong) parameter sequence on $R$.
\end{enumerate}
\end{prop}

{\it Proof:} By Proposition \ref{Wolmer}, $\bold x$ is weakly
proregular on $S$ if and only if $j(\bold x)$ is weakly proregular
on $R$.  As in the proof of Proposition \ref{Wolmer}, it suffices to
prove for sequences $\bold y$ of $R$ that $\bold y$ is a parameter
sequence on $S$ if and only if it is a parameter sequence on $R$.
Let $P$ be a prime minimal in $\operatorname{Supp}_S S/(\bold y)S$.
Since $S_{P}\cong R_{j(P)}\times (M_i)_{j(P)}$, we may assume
$\sqrt{(\bold y)R}=m$. It suffices to show that $H^\ell_{\bold
y}(R)\cong H^\ell_{\bold y}(S)$ where $\ell=\ell (\bold y)\ge 1$.
From the short exact sequence of $R$-modules $0\to M\to
S\xrightarrow{j} R\to 0$ we obtain the exact sequence
$$\cdots \to H^{\ell}_{\bold y}(M)\to H^{\ell}_{\bold y}(S)\to
H^{\ell}_{\bold y}(R)\to 0.$$ Since $H^k_{\bold y}(k(p))=0$ for
$k\ge 1$ and all $p\in \Spec R$, we have $H^{\ell}_{\bold
y}(M_i)=0$. The result now follows from the long exact sequence
above. \qed

\medskip

The following special case of Proposition \ref{sps-trivial} will be
needed in the proof of Example \ref{maximal}:

\begin{cor} \label{trivial-sps} Let $(R,m)$ be a Noetherian local ring and let $S=R\times
M_i$ as in Proposition \ref{Wolmer}.  Let $\bold x$ be a finite
sequence of elements of $S$.  Then $\bold x$ is a strong parameter
sequence on $S$ if and only if $\hgt (x_1,\dots,x_i)S=i$ for
$i=1,\dots,\ell(\bold x)$. \end{cor}

\section{The Cohen-Macaulay property}

We begin this section with a definition of Cohen-Macaulay for
arbitrary commutative rings:

\begin{df}{\rm A ring $R$ is called {\it Cohen-Macaulay} if every
strong parameter sequence on $R$ is a regular sequence.}
\end{df}

It follows from Remark \ref{Noeth-remark} that this definition
agrees with the usual definition of Cohen-Macaulay for Noetherian
rings (see \cite{BH} or \cite{Mat}).  Of course, given the many
different characterizations of Noetherian Cohen-Macaulayness, there
are many choices for extending the concept to non-Noetherian rings.
As we will see below, the property defined in Definition 4.1 has
many similarities to Noetherian Cohen-Macaulayness as well as some
stark differences (cf. Examples \ref{maximal} and
\ref{powerseries}). However, we believe the present definition is a
good one for exploring homological properties of rings, particularly
rings associated in some way to regular rings (e.g., invariant
subrings of regular rings).

Below we give equivalent formulations of Cohen-Macaulay in terms
of the polynomial grade, Koszul homology, and \u Cech cohomology
of strong parameter sequences:

\begin{prop} \label{CMGrade} Let $R$ be a ring.  The following
conditions are equivalent: \begin{enumerate}[(a)] \item $R$ is
Cohen-Macaulay. \item $\grade (\bold x)R=\ell (\bold x)$ for every
strong parameter sequence $\bold x$ of $R$. \item $\pgrade (\bold
x)R=\ell (\bold x)$ for every strong parameter sequence $\bold x$
of $R$. \item $H_i(\bold x;R)=0$ for all $i\ge 1$ for every strong
parameter sequence $\bold x$ of $R$.\item $H^i_{\bold x}(R)=0$ for
all $i<\ell(\bold x)$ for every strong parameter sequence $\bold
x$ of $R$.
\end{enumerate}
\end{prop}

{\it Proof:}  It is clear that
(a)$\Rightarrow$(b)$\Rightarrow$(c). By Proposition
\ref{CechGrade}, we also have that (c), (d), and (e) are
equivalent. It suffices to show (c) $\Rightarrow$ (a). We proceed
by induction on $\ell (\bold x)$ to show that all strong parameter
sequences on $R$ are regular sequences. If $\pgrade (x_1)R=1$ then
$x_1$ is a regular sequence. Suppose all strong parameter
sequences of length at most $\ell-1$ are regular sequences on $R$.
Let $\bold x$ be a strong parameter sequence of $R$ of length
$\ell$. Then by the induction hypothesis $\bold x'$ is a regular
sequence on $R$. Let $R'=R/(\bold x')R$. By Proposition
\ref{polygrade}(a), $\pgrade((x_{\ell})R',R')=1$, which implies
$x_{\ell}$ is a regular element on $R'$.  Hence, $\bold x$ is a
regular sequence on $R$. \qed

\medskip

The following example shows that, contrary to the Noetherian case,
it is not sufficient that $\pgrade (\bold x)R=\ell (\bold x)$ for
all {\it maximal} strong parameter sequences for a ring $R$ to be
Cohen-Macaulay

\begin{ex} \label{maximal} Let $(R,m)$ be a Cohen-Macaulay local ring of dimension $d>0$
and $S=R\times M_{d-1}$ as in Proposition 2.8.  Then $\pgrade (\bold
x)R=\ell (\bold x)$ for all maximal strong parameter sequences of
$S$, but $S$ is not Cohen-Macaulay.  In fact, $\pgrade (\bold x)R=0$
for all parameter sequences $\bold x$ of length less than $\dim R$.
\end{ex}

{\it Proof:} Combine Proposition \ref{Wolmer}, Corollary
\ref{trivial-sps}, and Proposition \ref{CMGrade} \qed

\medskip

We make some elementary observations concerning Cohen-Macaulay rings
of small dimension:

\begin{prop} Let $R$ be a ring.
\begin{enumerate}[(a)]
\item If $\dim R=0$ then $R$ is Cohen-Macaulay. \item If $R$ is a
one-dimensional domain, then $R$ is Cohen-Macaulay.
\end{enumerate}
\end{prop}

{\it Proof:} Part (a) follows from the fact that a
zero-dimensional ring has no parameter sequences.  For part (b),
note that by Proposition \ref{vanishing}, the maximum length of a
parameter sequence is one.  Since every non-zero element of $R$ is
a non-zero-divisor, we conclude that $R$ is Cohen-Macaulay. \qed

\medskip

The Cohen-Macaulay property descends along faithfully flat
extensions:

\begin{prop} Let $f:R\to S$ be a faithfully flat ring
homomorphism.  If $S$ is Cohen-Macaulay, then so is $R$.
\end{prop}

{\it Proof:}  This follows from Proposition \ref{param-prop}, part
(c) and \cite[Theorem 7.5]{Mat}. \qed

\medskip

One immediate consequence is:

\begin{cor} Let $R$ be a ring such that the polynomial ring $R[t]$
is Cohen-Macaulay.  Then $R$ is Cohen-Macaulay.
\end{cor}

We do not know whether $R[t]$ must be Cohen-Macaulay whenever $R$
is.  Likewise, we do not know whether the Cohen-Macaulay property
localizes.   In both cases, the difficulty lies in linking strong
parameter sequences of $R[t]$ (or $R_S$) to strong parameter
sequences of $R$. It may be that some mild condition on the ring,
such as requiring that the sets of minimal primes of finitely
generated ideals are finite, is necessary for these properties to
hold (cf. \cite{Ma}).  However, we do have the following:

\begin{prop} \label{locallyCM} Let $R$ be a ring and suppose $R_m$ is Cohen-Macaulay
for all maximal ideals $m$ of $R$.  Then $R$ is Cohen-Macaulay.
\end{prop}

{\it Proof:}  We proceed by induction on the length of a strong
parameter sequence $\bold x$ on $R$.  If $\ell(\bold x)=1$ then
$x_1$ is a regular sequence on $R_m$ for all maximal ideals $m$
containing $x_1$. Hence $x_1$ is regular on $R$.  Now assume
$\ell=\ell(\bold x)\ge 2$ and that all strong parameter sequences
on $R$ of length $\ell -1$ are regular sequences. For all maximal
ideals $m$ of $R$ containing $(\bold x)$, $x_{\ell}$ is regular on
$R_m/(\bold x')R_m$. Hence, $x_{\ell}$ is regular on $R/(\bold
x')R$. \qed

\medskip

We will call a ring $R$ {\it locally Cohen-Macaulay} if $R_p$ is
Cohen-Macaulay for all $p\in \Spec R$.  By Proposition 4.7, if $R$
is locally Cohen-Macaulay then $R_S$ is Cohen-Macaulay for all
multiplicatively closed sets $S$ of $R$.  As the following theorem
shows, coherent regular rings are locally Cohen-Macaulay.  A ring is
{\it regular} if every finitely generated ideal has finite
projective dimension (\cite{Be}). A ring is {\it coherent} if every
finitely generated ideal of $R$ is finitely presented.  (See
\cite{G1} for basic properties of coherent rings.) It is easily seen
that every finitely generated ideal of a coherent regular ring has
an FFR.

\medskip

\begin{thm} \label{regular} Let $R$ be a coherent regular ring.  Then $R$ is
locally Cohen-Macaulay.
\end{thm}

{\it Proof:} Since the localization of a coherent regular ring is
again coherent regular, it suffices to prove that $R$ is
Cohen-Macaulay.  Let $\bold x$ be a strong parameter sequence on
$R$ and $I=(\bold x)R$. By induction on $\ell=\ell(\bold x)$, we
may assume $\bold x'$ is a regular sequence on $R$.  We will
suppose $x_{\ell}$ is a zero-divisor on $R'=R/(\bold x')R$ and
derive a contradiction. Then $x_{\ell}\in p$ for some $p\in \Ass
R'$.  By Lemma \ref{attached} we have $\pdepth R'_p=0$.   Thus,
$\pdepth R_p =\ell-1$. By localizing at $p$, we can assume $(R,m)$
is a coherent regular quasi-local ring and $\pdepth R=\ell-1$. The
Auslander-Buchsbaum formula (Proposition \ref{ABT}) yields $\pd_R
R/I^t\le \pdepth R = \ell-1$ for all $t\ge 1$. Therefore,
$\Ext^{\ell}_R(R/I^t,R)=0$ for all $t\ge 1$.  Taking direct limits
we obtain $H^{\ell}_I(R)=0$.  By Proposition \ref{Schenzel}, we
have $H^{\ell}_{\bold x}(R)=0$, contradicting that $\bold x$ is a
parameter sequence on $R$. \qed

\medskip

We note that this theorem answers the question of Glaz  (\cite[p.
220]{G2}) mentioned in the introduction:  Does there exist a
definition of Cohen-Macaulay which agrees with the usual notion in
the Noetherian case and having the property every coherent regular
ring is Cohen-Macaulay?

The following is an example of a two-dimensional Cohen-Macaulay
quasi-local ring $R$ with the property that $R/xR$ is not
Cohen-Macaulay for some non-zero-divisor $x$ on $R$.

\begin{ex} \label{powerseries} {\rm Let $S=\mathbb C[[x,y]]$ be the
ring of formal power series in $x$ and $y$ over the field of
complex numbers.  Let $R=\mathbb C+ x\mathbb C[[x,y]]\subseteq S$.
It is easily seen that $R$ is a quasi-local domain.  We prove
that:
\begin{enumerate}
\item $R$ is Cohen-Macaulay. \item $R/xyR$ is not Cohen-Macaulay.
\end{enumerate}}
\end{ex}

{\it Proof:} Let $m=x\mathbb C[[x,y]]$ denote the maximal ideal of
$R$.  As $R$ is a domain every non-zero element of $m$ is both a
parameter and a non-zero-divisor.  To prove $R$ is Cohen-Macaulay,
it suffices to prove that $R$ has no parameter sequences of length
greater than one.   In fact, we will show $H^i_{\bold w}(R)=0$ for
all $i\ge 2$ and for all finite sequences $\bold w$ of $R$.  Let
$\bold w$ be such a sequence.  Clearly, we may assume $(\bold
w)\subseteq m$.   Consider the short exact sequence of $R$-modules
$$0\to R\to S\to S/R\to 0.$$
Note that as $xS\subseteq R$ we have $m(S/R)=0$.  Hence,
$H^i_{\bold w}(S/R)=0$ for all $i\ge 1$.  Therefore, $H^i_{\bold
w}(R)\cong H^i_{\bold w}(S)$ for all $i\ge 2$.   By the change of
rings theorem and since $S$ is Noetherian, we have that
$H^i_{\bold w}(S)\cong H^i_{(\bold w)S}(S)$ for all $i$.  Now,
$(\bold w)S\subseteq xS$ and so $\dim S/(\bold w)S>0$.  Since $S$
is a complete local domain of dimension two, we have $H^i_{(\bold
w)S}(S)=0$  for all $i\ge 2$ by the Hartshorne-Lichtenbaum
vanishing theorem \cite[Theorem 8.2.1]{BS}.  Hence, $H^i_{\bold
w}(R)=0$ for all $i\ge 2$.

Clearly, $xy$ is a non-zero-divisor on $R$.  We claim that $x$ is a
parameter on $R/xyR$.  Since $m=\sqrt{xR}$, we need only check that
$H^1_{x}(R/xyR)\neq 0$ and that $x$ is weakly proregular on $R/xyR$.
Clearly $x$ is not nilpotent in $R/xyR$ and so $H^1_{x}(R/xyR)\neq
0$ by Lemma \ref{easylemma}.  Also, $(xyR:_Rx)=xyS=(xyR:_Rx^2)$.  To
see this, first note that $xS\subseteq R$ and $yS\cap R=xyS$.  Thus,
$xyS$ is a prime ideal of $R$.  Since $xyR\subseteq xyS$ and
$x^2\not\in xyS$ we have that $(xyR:_Rx^2)\subseteq xyS$. On the
other hand, $x(xyS)=xy(xS)\subseteq xyR$ and thus $xyS\subseteq
(xyR:_Rx)\subseteq (xyR:_Rx^2)\subseteq xyS$.  Hence, $x$ is weakly
proregular on $R/xyR$.  As $x$ is a parameter and a zero-divisor on
$R/xyR$ (as $xy^2\in xyS\setminus xyR$), we see that $R/xyR$ is not
Cohen-Macaulay. \qed \medskip

We also note a connection between the present definition of
Cohen-Macaulay and the unmixedness notions proposed in \cite{Ha1}
and \cite{Ha2} as possible definitions of Cohen-Macaulay for
non-Noetherian rings.  An ideal $I$ of a ring $R$ is said to be \it
unmixed \rm if $\Ass R/I=\Min R/I$. It is well-known that a
Noetherian ring is Cohen-Macaulay if and only if every ideal
generated by a parameter sequence is unmixed \cite[Theorem
17.6]{Mat}. For arbitrary rings, this unmixedness condition implies
the notion of Cohen-Macaulay introduced here, but is properly
stronger.

\begin{prop} \label{unmixed}  Let $R$ be a ring such that every ideal generated by a
strong parameter sequence is unmixed.  Then $R$ is Cohen-Macaulay.
\end{prop}

{\it Proof:}  Let $\bold x$ be a strong parameter sequence on $R$.
We use induction on  $\ell(\bold x)$ to prove $\bold x$ is a
regular sequence. This is trivial in the case $\ell (\bold x)=0$.
Suppose $\ell=\ell(\bold x)>0$ and that $\bold x'$
 is a regular sequence on $R$.  It suffices to show that
 $x_{\ell}$ is regular on $R/(\bold x')R$.  If $x_{\ell}$ is a
 zero-divisor on $R/(\bold x')R$ then $x_{\ell}\in p$ for some
 $p\in \Ass R/(\bold x')R=\Min R/(\bold x')R$.  Thus,
 $\sqrt{(\bold x)R_p}=\sqrt{(\bold x')R_p}$.  Hence,
 $H^{\ell}_{\bold x}(R)_p=H^{\ell}_{\bold x'}(R)_p=0$,
 contradicting that $\bold x$ is a strong parameter sequence on
 $R$.  Thus, $\bold x$ is a regular sequence on $R$ and $R$ is
 Cohen-Macaulay. \qed \medskip

\medskip

Consequently, the weak Bourbaki unmixed rings and weak Bourbaki
height-unmixed rings studied in \cite{Ha1} and \cite{Ha2} are
Cohen-Macaulay.  However, if $R$ and $xy$ are as in Example
\ref{powerseries}, then $R$ is Cohen-Macaulay, $xy$ is a strong
parameter sequence on $R$, but $xyR$ is not unmixed (as $m\in \Ass
R/(xy)R\setminus \Min R/(xy)R$). Hence the converse of Proposition
\ref{unmixed} is false.

Let $R$ be an excellent Noetherian local domain of dimension $d$.
The {\it absolute integral closure} $R^+$ of $R$ is defined to be
the integral closure of $R$ in an algebraic closure of its field of
fractions. In \cite{HH2}, Hochster and Huneke prove that if
$\operatorname{char}R=p>0$ then $R^+$ is a big Cohen-Macaulay
algebra; i.e., every system of parameters for $R$ is a regular
sequence on $R^+$. Using this result, we can show that $R^+$ is a
Cohen-Macaulay ring in the sense introduced here.

\begin{thm} \label{Rplus} Let $R$ be an excellent Noetherian domain of
characteristic $p>0$.  Then $R^+$ is Cohen-Macaulay.
\end{thm}

{\it Proof:} Let $\bold x$ be a strong parameter sequence on
$R^+$. If we let $S=R[\bold x]$ then $S$ is also an excellent
Noetherian domain and $S^+=R^+$.  Hence, we may assume $\bold x$
is a sequence of elements in $R$. By Propositions \ref{CMGrade}
and \ref{CechGrade}, it suffices to prove that $H^i_{\bold
x}(R^+)=0$ for all $i<\ell(\bold x)$. Since integral closure and
\u Cech cohomology commute with localization, it suffices to prove
this in the case when $(R,m)$ is local and $(\bold x)R\subseteq
m$. By Proposition \ref{param-prop}(d), $\bold x$ is a strong
parameter sequence on $R$. Since $R$ is Noetherian, this means
$\bold x$ is a (partial) system of parameters for $R$.  By
\cite[Theorem 5.15]{HH2}, $\bold x$ is a regular sequence on
$R^+$. \qed

\medskip

As an application of non-Noetherian Cohen-Macaulayness, we consider
a conjecture raised by Glaz \cite{G2}:  {\it Let $R$ be a coherent
regular ring, $G$ a group of automorphisms of $R$, and $R^G$ the
ring of invariants.  Assume that there exists a module retraction
$\rho:R\to R^G$ and that $R$ is finitely generated $R^G$-module.
Then $R^G$ is Cohen-Macaulay.}  This conjecture is well-known to be
true in the case $R$ is Noetherian by the theorem of Hochster and
Eagon \cite[Proposition 12]{HE}. While we are not able to completely
resolve Glaz's conjecture using the present notion of
Cohen-Macaulay, we are able to prove it in the case $\dim R=2$
(Corollary \ref{invariant}).

Let $f:R\to S$ be a ring homomorphism.  A {\it module retraction}
from $S$ to $R$ is a $R$-module homomorphism  $\rho:S\to R$ such
that $\rho(f(r))=r$ for all $r\in R$.  In this case, we call $R$ a
{\it module retract} of $S$.

We begin with a basic lemma:

\begin{lem} \label{basic} Let $R\subseteq S$ be commutative rings such that $R$
is quasi-local, $S$ is finite over $R$, and there exists a module
retraction $\rho:S\to R$. Then there exists a maximal ideal $q$ of
$S$ such that $\rho(xS)=R$ for every $x\in S\setminus q$.
\end{lem}

{\it Proof:} Let $m$ be the maximal ideal of $R$. Since
$\rho(mS)\subseteq m$, $\overline{\rho}:S/mS\to R/m$ is a retraction
of the extension $R/m\subseteq S/mS$.  Thus, it suffices to prove
the lemma in the case $R=k$ is a field.  Let $q_1\cap \cdots \cap
q_t=0$ be the primary decomposition for $S$ where $p_i=\sqrt{q_i}$
for each $i$.  (Since $k$ is a field each $p_i$ is a maximal ideal
of $S$.) The Chinese Remainder Theorem gives an isomorphism $S\cong
\prod_i Se_i$ where $S/q_i\cong Se_i\subseteq S$ for $i=1,\dots,t$.
Since $1=\rho(1)=\rho(e_1)+\cdots +\rho(e_t)$, we have
$r_j=\rho(e_j)\neq 0$ for some $j$. Let $\phi_j:S/q_j\to S$ be
defined by $\phi_j(\overline{s})=se_j$. Setting
$\rho_j=r_j^{-1}\rho\, \phi_j$, we obtain a module retraction for
the extension $k\subseteq S/q_j$. Suppose the lemma holds for the
retraction $\rho_j$.  Then for each $x\in S\setminus p_j$ we have
$\rho_j(x(S/q_j))=k$. Thus $\rho(xS)\supseteq
\rho(xe_jS)=r_j\rho_j(x(S/q_j))=k$.  This reduces the proof of the
lemma to the case where $R$ is a field and $S$ is local.  But this
case is trivial, since if $x$ is a unit in $S$ then
$\rho(Sx)=\rho(S)=R$. \qed

\begin{cor} \label{faithful}  Let $R$, $S$, $\rho$, and $q$ be as in
Lemma \ref{basic}.  Let $M$ be an $R$-module and $\bold x$ a
sequence of elements from $R$.  The following hold:
\begin{enumerate}[(a)]
\item $M=0$ if and only if $S_q\otimes_RM=0$. \item If $H^i_{\bold x}(S\otimes_RM)_q=0$
then $H^i_{\bold x}(M)=0$.
\item For any ideal $I$ of $R$,
$$\pgrade (IS_q,S_q\otimes_RM)\le \pgrade(I,M).$$ \item If $\bold
x$ is a regular sequence on $S_q$ then $\bold x$ is a regular
sequence on $R$.
\end{enumerate}
\end{cor}

{\it Proof:} Observe that the map $\rho':S\otimes_RM\to M$ given
by $\rho'(s\otimes m)=\rho(s)m$ is a retraction of the map $j:M\to
S\otimes_RM$ given by $j(m)=1\otimes m$ for all $m\in M$. Suppose
$S_q\otimes_RM=0$ and let $m\in M$. Then there exists $t\in
S\setminus q$ such that $t\otimes m=t(1\otimes m)=0$.  By the
lemma, there exists $s\in S$ such that $\rho(st)=1$.  Then
$$m=\rho(st)m=\rho'(st\otimes m)=\rho'(0)=0.$$ To prove (b) we
first note that the maps $j$ and $\rho$ induce maps of complexes
$$C(\bold x;M)\xrightarrow{\psi} S\otimes_R C(\bold
x;M)\xrightarrow{\phi} R\otimes_RC(\bold x;M)\cong C(\bold x;M)$$
where $\psi$ is the obvious inclusion induced by $j$ and
$\phi=\rho\otimes 1$. This retraction of complexes induces maps on
homology $\phi^*_i:H^i_{\bold x}(S\otimes_RM)\to H^i_{\bold x}(M)$
where $\phi^*_i(su)=\rho(s)u$ for every $u\in H^i_{\bold x}(M)$.
(The injection $\psi$ allows us to view $H^i_{\bold x}(M)$ as a
direct summand of $H^i_{\bold x}(S\otimes_RM)$.)  Part (b) now
follows by the same argument as in the proof of part (a) with
$\phi^*_i$ in place of $\rho'$. For part (c), we note that we may
assume that $I=(\bold x)R$ is a finitely generated ideal by
Proposition \ref{polygrade}(c). The inequality is now immediate
from part (b) and Proposition \ref{CechGrade}. Part (d) follows
from part (c) and induction on the length of the regular sequence
as in the proof of Proposition \ref{CMGrade}. \qed

\medskip

We now apply these results to parameter sequences on $R$:

\begin{thm} \label{regular-retract} Let $R\subseteq S$ be
commutative rings such that $R$ is a module retract of $S$, $S$ is
finite over $R$, and $S$ is a coherent regular ring.  Let $\bold
x$ be a strong parameter sequence on $R$ such that $\ell(\bold
x)\le 2$. Then $\bold x$ is a regular sequence on $R$.
\end{thm}

{\it Proof:} Let $x$ be a parameter on $R$ and suppose $x$ is a
zero-divisor on $R$.  Then $x\in p$ for some $p\in \Ass R$.  Since
the image of $x$ is nonzero in $R_p$ (as $x$ is a parameter), we can
localize $R$ and $S$ at $p$ and assume $(R,m)$ is quasi-local. Let
$q$ be the maximal ideal of $S$ given by Lemma \ref{basic}. Then the
image of $x$ in $S_q$ is nonzero by part (a) of Proposition
\ref{faithful}. As $S_q$ is a coherent regular quasi-local ring, it
is a domain (\cite{Be}). Hence $x$ is a non-zero-divisor on $S_q$
and thus on $R$ as well by Proposition \ref{faithful}(c).

Now suppose $x,y$ is a strong parameter sequence on $R$.  By above,
we have that $x$ is a non-zero-divisor on $R$.  Suppose $y$ is a
zero-divisor on $R/xR$.  Then $y\in p$ for some $p\in \Ass R/xR$.
Localizing at $p$, we can assume that $(R,m)$ is quasi-local and
$\pdepth R/xR=0$.  Again, let $q$ be the maximal ideal of $S$ given
by Lemma \ref{basic}.  By the argument given in the parameter case,
$x$ is a non-zero-divisor on $S_q$.  Since $qS_q=\sqrt{mS_q}$,
\begin{align}
\pdepth S_q/xS_q&=\pgrade (mS_q,S_q/xS_q) \notag \\
&\le \pgrade (m,R/xR) \qquad \text{(by Corollary \ref{faithful}(b))}
\notag
\\&= 0. \notag
\end{align}
Hence, $\pdepth S_q=1$.  Therefore, $\pd_{S_q}S_q/(x,y)S_q\le 1$,
which implies that $(x,y)S_q$ is principal.  Thus,
$H^2_{x,y}(S)_q=0$ and hence $H^2_{x,y}(R)=0$ by Proposition
\ref{faithful}(b), contradicting that $x,y$ is a strong parameter
sequence.\qed

\medskip

As a special case, we get the following:

\begin{cor} \label{invariant} Let $R$ be a coherent regular ring
of dimension at most two and $G$ a finite group of automorphisms
of $R$ such that the order of $G$ is a unit in $R$. Let $R^G$ be
the subring of invariants of $R$ under the action of $G$ and
assume that $R$ is a finite $R^G$-module.  Then $R^G$ is a
coherent locally Cohen-Macaulay ring.
\end{cor}

{\it Proof:}  The averaging map gives a module retraction from $R$
to $R^G$.  Furthermore, $(R^G)_p$ is a module retract of $R_p$ for
every prime $p$ of $R^G$.  Since $\dim (R^G)_p=\dim R_p\le 2$, the
maximal length of any parameter sequence on $(R^G)_p$ is two by
Proposition \ref{vanishing}.  Applying Theorem \ref{regular-retract}
we see that $(R^G)_p$ is Cohen-Macaulay. Coherence of $R^G$ follows
from \cite[Theorem 1]{G3}.  \qed

\end{document}